\documentclass[a4paper,11pt]{article}
\usepackage{graphicx}
\usepackage{amssymb}
\usepackage{amsmath}
\usepackage{mathtools}
\usepackage{amsfonts}
\usepackage{appendix}
\usepackage{color}
\usepackage{multirow}
\usepackage{pdflscape}
\usepackage{caption}
\usepackage{subcaption}
\usepackage{xcolor}
\usepackage{setspace}
\usepackage{booktabs}
\usepackage{adjustbox}
\usepackage{float}
\def\cov{\mathop{\rm Cov}\nolimits}

\begin{document}
	
\begin{center}
{\Large \bf A sure independence screening procedure for ultra-high dimensional partially linear additive models }
\end{center}

\begin{center}
{\bf M. Kazemi, D. Shahsavani, M. Arashi }

\vspace*{0.4 cm}

{\it  Department of Statistics, Shahrood University of Technology, Shahrood, Iran.}
\end{center}

\bigskip
\begin{abstract}\
 We introduce a two-step procedure, in the context of ultra-high dimensional additive models, which aims to reduce the size of  covariates vector and distinguish linear and nonlinear effects among nonzero components. Our proposed screening procedure, in the first step, is constructed based on the concept of cumulative distribution function and conditional expectation of response in the framework of marginal correlation. B-splines and empirical distribution functions are used to estimate the two above measures. The sure property of this procedure is also established.  In the second step, a double penalization based procedure is applied to identify nonzero and linear components, simultaneously. The performance of the designed method is examined by several test functions to show its capabilities against competitor methods when errors distribution are varied. 
Simulation studies imply that the proposed screening procedure can be applied to the ultra-high dimensional data and well detect the influential covariates.  It is also demonstrate the superiority in comparison with the existing methods. This method is also applied to identify most influential genes for overexpression of a G protein-coupled receptor in mice.


\end{abstract}

\section{Introduction}
Dimension reduction and identifying the relevant vector components are challenges in prediction problems. Many endeavours have been made to identify the irrelevant components via screening or variable selection methods. A group of statisticians considered this subject in the framework of the partially linear additive model with form
\begin{equation}\label{q1}
Y=\sum_{j\epsilon S_1}\beta_j X_{j}+\sum_{j \epsilon S_2}f_j(X_{j})+\varepsilon,
\end{equation}
for data pairs $ (y_i,x_{i1},\dots x_{ip}),1\leq i \leq n $, where $ Y $ is the response and a $ p $-dimensional covariate vector $ X=( X_{1},\dots,X_{p}) $  is divided into two mutually exclusive and complementary subsets $ S_1 $ and $ S_2 $.
It is also assumed that the mean response is linearly related to the covariates in $ S_1 $, with the regression coefficients of $ \{\beta_j:j\epsilon S_1\} $, and the the remaining covariates in $ S_2 $ are included in the nonparametric part of the model through smooth functions $ \{f_j:j\epsilon S_2\} $ and the model error $ \varepsilon $ has conditional mean zero and finite variance $ \sigma^2 $ given $ X $. To ensure identifiability of the nonparametric functions, it is assumed usually that $  E[f_j(X_j)] = 0 $ for $ j\epsilon S_2 $. Estimation and variable selection for partially linear additive models have been well studied in literature, and we refer, for instance, to Liu et al. (2011), Lian (2012a), Guo et al. (2013), Du et al. (2015), Lv et al. (2016), among others.\

The use of model (\ref{q1}) is based on the assumption that the linear and nonlinear parts are known in advance. However, such prior information is usually unavailable, especially when the number of covariates is large. Thus, in addition to distinguish nonzero components, it is of great interest to develop some efficient methods to identify linear components from nonlinear ones. For this reason, our attention in this article is focused on general additive models
\begin{equation}\label{q2}
Y=\sum_{j=1}^{p}f_j(X_j)+\epsilon.
\end{equation}

Zhang et al. (2011) studied the model selection using two penalties, simultaneously, to identify the zero and linear components in partially linear additive models. Their method is selection consistent in the special case of tensor product design. However, they did not prove any selection consistency results for general partially linear models. Motivated by this, Huang et al (2012) proposed a semiparametric regression pursuit method for distinguishing linear from nonlinear components using a group MCP penalty and showed that the proposed approach is model-pursuit consistent. Lian (2012b) provided a way to determine linear components by using SCAD penalty using B-spline expansion. This was a new usage of SCAD in which no variable selection is performed. Lian (2012c) successfully identified nonzero and linear components of model (\ref{q2}) by applying a two-fold SCAD penalty in the additive quantile regression. \

When the number of covariates can diverge with the sample size, another two penalty procedure in high dimensional setting was proposed by Lian et al. (2015) in which insignificant predictors and parametric components were simultaneously identified in additive models. We will be using their method in our paper, however we have a different concern.\

Nevertheless, when $ p $ grows exponentially with $ n $, the aforementioned penalized variable selection methods may not work for the ultra-high dimensional partially linear additive model (\ref{q1}) due to the simultaneous challenges of computational expediency, statistical accuracy and algorithm stability (Fan et al., 2009). To address these challenges, sure independence screening (SIS) was introduced by Fan and Lv (2008) in the context of linear regression models for feature screening in ultrahigh-dimensional data analysis. Many authors further developed the SIS method and applied it to various statistical models, such as generalized linear models (Fan et al., 2009; Fan and Song, 2010), nonparametric additive models (NIS, Fan et al., 2011) and varying coefficient models (Fan et al., 2014; Liu et al., 2014). Furthermore, in order to avoid the specification of a particular model structure, Zhu et al. (2011) proposed a sure independent ranking and screening (SIRS) procedure for ultrahigh-dimensional data in the framework of the general multi-index models. Thereafter, a model-free SIS based on the distance correlation was developed by Li et al. (2012a). Using the Kendall $ \tau $,  Li et al. (2012b) proposed a robust screening procedure in the framework of the transformation models. In a model-free fashion, Zhang et al. (2017) proposed a correlation rank screening procedure (CR-SIS), which can naturally handle ultrahigh-dimensional survival data based on the covariance between unconditional distribution function of $ Y $ and covariates.\

In this paper, a screening procedure is followed by structure identification and variable selection method. We apply the modified version of Zhang et al. (2017) to reduce the dimensionality in ultra-high dimensional partial linear additive models, and then use the double penalization based procedure of Lian et al. (2015) to simultaneously identify nonzero and linear components.\

The plan of paper is as follows. In Section 2, a modification of nonparametric independence screening procedure Zhang et al. (2017) is introduced and its theoretical properties are considered. After performing dimension reduction, the doubly penalized estimation method of Lian et al. (2015) is explained in details in Section 3. In section 4, simulation studies are carried out to assess the performance of the proposed method and to compare it with some existing methods. A real data example is used for illustration in Section 5.

\section{Screening Procedure}

In the category of model-free screening procedures for ultrahigh dimensional setting, a correlation based sure independence screening method (CR-SIS) was suggested by Zhang et.al (2017) from distribution function prospect to reduce the dimension of potential covariates. Assume $ E(X_j)=0 $ and define the active covariate set as
\begin{equation}
\mathcal{A}=\{1\leq j\leq p: F(y|X) ~\text{depends on} ~ X_j\},
\end{equation}
to identify the contribution of each covariate to the distribution function of $ Y $ given $ X $, i.e., $ F(y|X)=P(Y\leq y|). $ To this aim, they considered the covariance between each covariate $ X_j $ and the unconditional distribution function of $ Y $, i.e.,
$ R_{j}(Y) = \cov\big(X_j, G(Y)\big)=E\big(X_jG(Y)\big) $
where  $ G(y)=P(Y\leq y) $ contains the whole information of $ Y $. Therefore the relationship between $ Y $ and $ X_j  $ could be reflected by the population version of marginal utility measure
\begin{equation*}
r_j = [R_j(Y)]^2,~j=1,\dots,p.
\end{equation*}

An estimator of $ r_j $ based on the random sample
$ (X_i,Y_i )_{i=1}^n , $ is given by
\begin{equation}
\hat{r}_j=\Big\{\frac{1}{n}\sum_{i=1}^n X_{ij} \hat{G}_n (Y_i)\Big\}^2
\end{equation}
where
\begin{equation*}
\hat{G}_n(y)=\frac{1}{n}\sum_{i=1}^{n}I(Y_i\leq y)
\end{equation*}
is the empirical distribution function. Thus, by ranking the $ r_k $ from largest to smallest, the important predictors are determined by the estimated active set,
\begin{equation}
\hat{\mathcal{A}}=\left\lbrace 1\leq j\leq p:  \hat{r}_j \geq cn^{-\alpha} \right\rbrace,
\end{equation}
for some constants $ c > 0 $ and $ \alpha \in [0, 1/2) $.\

Our modification of the above method is based on the approach of Fan et al. (2011) for the additive model (\ref{q2}). \

Consider $ p $ marginal nonparametric regression models
$  Y=f_j(X_j)+\epsilon,j=1,\dots ,p $
and obtain $ m_j=E(Y|X_j) $ as the solution of minimization problems:
\begin{equation*} \label{q4}
\min_{m_j\epsilon L_2(P)} E\big[Y-m_j(X_j)\big]^2,
\qquad j=1\dots p,
\end{equation*}
where $ L_2(P) $ is the class of square integrable functions under the measure $ P $. Then, the measure $ E(m_j (X_j))^2 $ is used for ranking the utility of covariates in the model (\ref{q2}). It can be noted that
\begin{equation*}
\cov(Y,m_j (X_j))=E(Ym_j (X_j))=E(m_j (X_j))^2.
\end{equation*}

The key measure $ R_j (Y)=\cov(X_j,G(Y)) $ in the model free fashion can be modified by considering the nonparametric measure $ \cov(m_j (X_j ),Y) $ in additive models. Our new measure of correlation between $ Y $ and $ X_j $ is proposed by substituting
$ m_j $ instead of $ X_j $ in the $ R_j $ formula:
\begin{equation}
\cov(m_j,G(Y))=\cov(E(Y|X_j ),G(Y)).
\end{equation}

The estimation of $ m_j $ can also be done by using the B-spline functions. Let $ \{B_{j1}(x),...,B_{jK}(x)\} $  be the normalized B-spline basis functions of order $ q $. With this, we have the following approximation
\begin{equation*}
\hat{m}_{j}(x)=B_{j}(x)^{T}\hat{\beta_{j}}, \quad 1\leq j\leq p,
\end{equation*}
where $ B_{j}(x) = \big(B_{j1}(x),...,B_{jK}(x)\big)^T $ and $\hat{\beta}_j=\left(\hat{\beta}_{j1},...,\hat{\beta}_{jK}\right)^T  $ is obtained through the componentwise least squares regression:
\begin{equation*}
\hat{\beta}_{j}= \text{arg}\min_{\beta_j\epsilon \mathbb{R}^{K} }\frac{1}{n} \sum_{i=1}^{n}\left[y_i-\beta_{j}^{T}B_{j}(x_{ij})\right]^2.
\end{equation*}
This procedure reduces the dimensionality from $  p $ to a possibly much smaller space with model size $ d=|\hat{\mathcal{A}}| $. The question is whether the procedure has a sure screening property, as postulated by Fan and Lv (2008).\

Here, we show that the proposed screening procedure possesses sure screening property. We impose the following regularity conditions throughout our discussion.\\

C1. There exists a positive constant $ \xi $ such that
\begin{equation*}
\max_{1\leq k\leq p} E[m_j(X_j)]^2<\xi.
\end{equation*}

C2. It holds that
\begin{equation*}
\min_{k\in \mathcal{A}}r_k\geq2cn^{-\alpha},
\end{equation*}

for some constants $ c > 0 $ and $ \alpha \in [0, 1/2). $\\

\textbf{Theorem 1}.
Under condition C1, there exists a constant $  \eta > 0 $ such that
\begin{equation*}
p(\max_{1\leq k\leq p}|\hat{r_k}-r_k|\geq cn^{-\alpha})\leq
O\big[p\exp{-\eta (\frac{n^{1-2\alpha}}{\log\log n})^{\frac{1}{2}}} \big].
\end{equation*}
Under conditions C1, C2, it holds that
\begin{equation*}
p(\mathcal{A}\subseteq \hat{\mathcal{A}})\geqslant 1-
O\big[a_n\exp{-\eta (\frac{n^{1-2\alpha}}{\log\log n})^{\frac{1}{2}}} \big],
\end{equation*}
where $ an = |\mathcal{A}| $ is the cardinality of $ \mathcal{A} $.\\
 
\section{Group Penalization}
Screening is an efficient method to reduce the model size from a very large value $ p $ to a moderate scale $ d $
by specifying sensible threshold parameters $ \nu_n $, whereas it is difficult to choose in practice. A practical way is to select the top $ d $ variables by ranking marginal utilities. The choice of $ d $ plays a very important role in the screening stage. Fan and Lv (2008) recommended $ d = [⌊n/log (n)] $ as a sensible choice. Such a $ d $ value is also suggested by Fan et al. (2009), which showed that the model-based, rather than data-driven, choice of $ d $ provides satisfactory and robust performance. Zhao and Li (2012) proposed an approach to select $ d $ for Cox models by controlling false positive rate. In this study, we adopt Fan et al. (2009)'s recommendation. A larger value of the specified $ d $ would give a greater chance to include inactive variables. This can be solved by a penalty-based variable selection procedure given below.\


Now, suppose that $ d $ variables are selected in the screening stage. Consider a joint nonparametric additive model
$
Y=\sum_{j=1}^{d}f_j(X_j)+\epsilon.
$
As in Section 2, B-spline basis is used to approximate each of unknown smooth functions, i.e., $ f_{j}(x)\approx\sum_{k}b_{jk}B_{jk}(x) $ for $ j=1,...,d $. We use the two-fold penalization procedure of Lian et al. (2015) to automatically identify different types of components, i.e., we find coefficient $ b=(b_{1}^{T},\dots b_{d}^{T})^T, b_j=(b_{j1},\dots b_{jK})^T, j=1, \dots, d,$
\begin{align}\label{q5}
\hat{b}=\text{arg min}_{b} \frac{1}{2}\sum_{i=1}^{n}
\left(  Y_i-\mu-\sum_{j=1}^{d}\sum_{k=1}^{K}b_{jk}B_{jk}(X_{ij})\right) ^{2} \nonumber \\
+n\sum_{j=1}^{d}p_{\lambda_1}(w_{1j} \rVert b_j \rVert_{A_j}) +n\sum_{j=1}^{d}p_{\lambda_2}(w_{2j} \rVert b_j \rVert_{D_j}),
\end{align}
where $ p_{\lambda}(|t|)=\lambda |t| $ is the LASSO penalty function, $ \lambda_1, \lambda_2 $ are regularization parameters. $ A_j $ and $ D_j $ are two $ K\times K $ matrices, $  \rVert b_j \rVert_{A_j}=\left(b_{j}^T A_{j}b_{j} \right)^{\frac{1}{2}},\rVert b_j \rVert_{D_j}=\left(b_{j}^T D_{j}b_{j} \right)^{\frac{1}{2}} $.
There is some flexibility in choosing  $ A_j $ and $ D_j $ but one requirement is that
$ \rVert b_j \rVert_{A_j}=0 $ if only if $\sum_{k}b_{jk}B_{jk}(x)\equiv 0 $
and
$ \rVert b_j \rVert_{D_j}=0 $ if only if $\sum_{k}b_{jk}B_{jk}(x)\equiv 0 $
is a linear function, so that the two penalties can be used to
identify zero and linear components, respectively. One natural choice is $ A_{j}= \{\int_{0 }^{1}B_{jk}(x)B_{jk'}(x) dx \}_{k,k'=1}^{K}$ and $ D_{j}= \{\int_{0 }^{1}B''_{jk}(x)B''_{jk'}(x) dx \}_{k,k'=1}^{K}$
so that
$ \rVert b_j \rVert_{A_j}= \rVert \sum_{k}b_{jk}B_{jk}(x) \rVert $
and
$ \rVert b_j \rVert_{D_j}= \rVert \sum_{k}b_{jk}B''_{jk}(x) \rVert $.
The adaptive group lasso penalty in (\ref{q5}) involves the weights vectors $ w_1 = (w_{11},...,w_{1d})$ and $ w_2 =(w_{21},...,w_{2d})$.
The weights $ w_{1j} $ are best if large for zero components and small for nonzero
ones, and similarly the  $ w_{2j} $ are best if large for linear components and small for
nonparametric ones. Using the group Lasso of Huang et al. (2010), the initial estimator is obtained as

\begin{equation}\label{q8}
\tilde{b}=\text{arg min}_{b} \frac{1}{2} \rVert Y-Zb \rVert ^2
+n \lambda_0\sum_{j=1}^{d}\rVert b_j \rVert_{A_j}
\end{equation}

Using this initial estimator, we can then set $ w_{1j}=\frac{1}{\rVert \tilde{b}_j \rVert_{A_j}} $
and
$ w_{2j}=\frac{1}{\rVert \tilde{b}_j \rVert_{D_j}} $ in (\ref{q5}). Let
$$
Z_j=\left(\begin{array}{cccc}
B_{j1}(X_{1j}) &  B_{j2}(X_{1j}) & \cdots & B_{jK}(X_{1j}) \\
\vdots         &        \vdots   &        & \vdots \\
B_{j1}(X_{nj}) & B_{j2}(X_{nj})  & \cdots & B_{jK}(X_{nj})
\end{array}  \right)_{n \times K},
$$
$ Z=(Z_1,...,Z_d) $ and $ Y=(Y_1,...,Y_n)$. Then (\ref{q5}) can be written in matrix form as
\begin{equation}\label{q6}
\hat{b}=\text{arg} \min_{b} \frac{1}{2} \rVert Y-Zb \rVert ^2
+n\sum_{j=1}^{d}p_{\lambda_1}(w_{1j} \rVert b_j \rVert_{A_j}) +n\sum_{j=1}^{d}p_{\lambda_2}(w_{2j} \rVert b_j \rVert_{D_j}).
\end{equation}
To find the minimum of (\ref{q6}) for fixed tuning parameters, we use the iterative local quadratic approximation (LQA) proposed by Fan and Li (2001). Using a simple Taylor expansion, given an initial estimate $ b_{j}^{0} $, if $ \rVert b_j \rVert_{A_j}>0 ~~ \text{and}~~ \rVert b_j \rVert_{D_j}>0 $,
we approximate the penalty terms by
\begin{equation*}
p_{\lambda_1}(w_{1j} \rVert b_j \rVert_{A_j})\approx p_{\lambda_1}\left(w_{1j} \rVert b_{j}^{(0)} \rVert_{A_j}\right)+
\frac{1}{2}\frac{p'_{\lambda_1}\left( w_{1j}\rVert b_{j}^{(0)} \rVert_{A_j}\right)}{w_{1j} \rVert b_{j}^{(0)} \rVert_{A_j}}
\left\lbrace w_{1j}^{2}\rVert b_{j} \rVert_{A_j}^{2}-w_{1j}^{2}\rVert b_{j}^{(0)} \rVert_{A_j}^{2}\right\rbrace,
\end{equation*}
and
\begin{equation*}
p_{\lambda_2}\left(w_{2j} \rVert b_j \rVert_{D_j}\right)\approx p_{\lambda_2}\left(w_{2j} \rVert b_{j}^{(0)} \rVert_{D_j}\right)+
\frac{1}{2}\frac{p'_{\lambda_2}\left( w_{2j}\rVert b_{j}^{(0)} \rVert_{D_j}\right)}{w_{2j} \rVert b_{j}^{(0)} \rVert_{D_j}}
\left\lbrace w_{2j}^{2}\rVert b_{j} \rVert_{D_j}^{2}-w_{2j}^{2}\rVert b_{j}^{(0)} \rVert_{D_j}^{2}\right\rbrace.
\end{equation*}
After removing some irrelevant terms, the criterion becomes
\begin{equation}\label{q7}
Q(b)=\frac{1}{n} \rVert Y-Zb \rVert ^2 +\frac{1}{2}b^{T}(\Omega_1+\Omega_2)b
\end{equation}
for two $ dK \times dK $ matrices $ \Omega_1 $ and $ \Omega_2 $ defined by
\begin{equation}
\Omega_1=\text{diag}\left( \frac{\lambda_{1}w_{11}}{\rVert b_{1}^{(0)} \rVert_{A_1}}A_1,\dots,
\frac{\lambda_{1}w_{1d}}{\rVert b_{d}^{(0)} \rVert_{A_d}}A_d \right) \nonumber
\end{equation}
and
\begin{equation}
\Omega_2=\text{diag}\left( \frac{\lambda_{2}w_{21}}{\rVert b_{1}^{(0)} \rVert_{D_1}}D_1,\dots,
\frac{\lambda_{1}w_{2d}}{\rVert b_{d}^{(0)} \rVert_{D_d}}D_d \right) \nonumber
\end{equation}

Note that (\ref{q7}) is a quadratic function and thus there exists a
closed-form solution. Then the updating equation given the current estimate
$ b^{(0)} $ is
\begin{equation}\label{key}
b=\big(Z^{T}Z+n(\Omega_1+\Omega_2 )\big)^{-1}Z^{T}Y
\end{equation}

The algorithm repeatedly solves the minimization criterion (\ref{q7}) and updates $ b^{(m)} $ to $ b^{(m+1)} $, $ m=0,1,... $ until convergence. That is, in the m-th iteration, we solve (\ref{q7}), where $ \Omega_1 $ and $ \Omega_2 $ are as defined above but with $ b_{j}^{0} $ replaced by the current estimate $ b_{j}^{(m)} $. The solution obtained from (\ref{q7}) is the new estimate $ b^{(m+1)} $. \
During the iterations, as soon as some $ \rVert b_j \rVert_{A_j} $
(respectively, $ \rVert b_j \rVert_{D_j} $) drops below a certain threshold ($ 10^{-6} $ in our implementation), the component is identified as a zero function (respectively, linear function). Screening is an efficient method to reduce the model size from a very large value $ p $ to a moderate scale $ d $
by specifying sensible threshold parameters $ \nu_n $, whereas it is difficult to choose in practice. A practical way is to select the top $ d $ variables by ranking marginal utilities. The choice of $ d $ plays a very important role in the screening stage. Fan and Lv (2008) recommended $ d = [⌊n/log (n)] $ as a sensible choice. Such a $ d $ value is also suggested by Fan et al. (2009), which showed that the model-based, rather than data-driven, choice of $ d $ provides satisfactory and robust performance. Zhao and Li (2012) proposed an approach to select $ d $ for Cox models by controlling false positive rate. In this study, we adopt Fan et al. (2009)'s recommendation. A larger value of the specified $ d $ would give a greater chance to include inactive variables. This can be solved by a penalty-based variable selection procedure given below.\


Now, suppose that $ d $ variables are selected in the screening stage. Consider a joint nonparametric additive model
$
Y=\sum_{j=1}^{d}f_j(X_j)+\epsilon.
$
As in Section 2, B-spline basis is used to approximate each of unknown smooth functions, i.e., $ f_{j}(x)\approx\sum_{k}b_{jk}B_{jk}(x) $ for $ j=1,...,d $. We use the two-fold penalization procedure of Lian et al. (2015) to automatically identify different types of components, i.e., we find coefficient $ b=(b_{1}^{T},\dots b_{d}^{T})^T, b_j=(b_{j1},\dots b_{jK})^T, j=1, \dots, d,$
\begin{align}\label{q5}
\hat{b}=\text{arg min}_{b} \frac{1}{2}\sum_{i=1}^{n}
\left(  Y_i-\mu-\sum_{j=1}^{d}\sum_{k=1}^{K}b_{jk}B_{jk}(X_{ij})\right) ^{2} \nonumber \\
+n\sum_{j=1}^{d}p_{\lambda_1}(w_{1j} \rVert b_j \rVert_{A_j}) +n\sum_{j=1}^{d}p_{\lambda_2}(w_{2j} \rVert b_j \rVert_{D_j}),
\end{align}
where $ p_{\lambda}(|t|)=\lambda |t| $ is the LASSO penalty function, $ \lambda_1, \lambda_2 $ are regularization parameters. $ A_j $ and $ D_j $ are two $ K\times K $ matrices, $  \rVert b_j \rVert_{A_j}=\left(b_{j}^T A_{j}b_{j} \right)^{\frac{1}{2}},\rVert b_j \rVert_{D_j}=\left(b_{j}^T D_{j}b_{j} \right)^{\frac{1}{2}} $.
There is some flexibility in choosing  $ A_j $ and $ D_j $ but one requirement is that
$ \rVert b_j \rVert_{A_j}=0 $ if only if $\sum_{k}b_{jk}B_{jk}(x)\equiv 0 $
and
$ \rVert b_j \rVert_{D_j}=0 $ if only if $\sum_{k}b_{jk}B_{jk}(x)\equiv 0 $
is a linear function, so that the two penalties can be used to
identify zero and linear components, respectively. One natural choice is $ A_{j}= \{\int_{0 }^{1}B_{jk}(x)B_{jk'}(x) dx \}_{k,k'=1}^{K}$ and $ D_{j}= \{\int_{0 }^{1}B''_{jk}(x)B''_{jk'}(x) dx \}_{k,k'=1}^{K}$
so that
$ \rVert b_j \rVert_{A_j}= \rVert \sum_{k}b_{jk}B_{jk}(x) \rVert $
and
$ \rVert b_j \rVert_{D_j}= \rVert \sum_{k}b_{jk}B''_{jk}(x) \rVert $.
The adaptive group lasso penalty in (\ref{q5}) involves the weights vectors $ w_1 = (w_{11},...,w_{1d})$ and $ w_2 =(w_{21},...,w_{2d})$.
The weights $ w_{1j} $ are best if large for zero components and small for nonzero
ones, and similarly the  $ w_{2j} $ are best if large for linear components and small for
nonparametric ones. Using the group Lasso of Huang et al. (2010), the initial estimator is obtained as

\begin{equation}\label{q8}
\tilde{b}=\text{arg min}_{b} \frac{1}{2} \rVert Y-Zb \rVert ^2
+n \lambda_0\sum_{j=1}^{d}\rVert b_j \rVert_{A_j}
\end{equation}

Using this initial estimator, we can then set $ w_{1j}=\frac{1}{\rVert \tilde{b}_j \rVert_{A_j}} $
and
$ w_{2j}=\frac{1}{\rVert \tilde{b}_j \rVert_{D_j}} $ in (\ref{q5}). Let
$$
Z_j=\left(\begin{array}{cccc}
B_{j1}(X_{1j}) &  B_{j2}(X_{1j}) & \cdots & B_{jK}(X_{1j}) \\
\vdots         &        \vdots   &        & \vdots \\
B_{j1}(X_{nj}) & B_{j2}(X_{nj})  & \cdots & B_{jK}(X_{nj})
\end{array}  \right)_{n \times K},
$$
$ Z=(Z_1,...,Z_d) $ and $ Y=(Y_1,...,Y_n)$. Then (\ref{q5}) can be written in matrix form as
\begin{equation}\label{q6}
\hat{b}=\text{arg} \min_{b} \frac{1}{2} \rVert Y-Zb \rVert ^2
+n\sum_{j=1}^{d}p_{\lambda_1}(w_{1j} \rVert b_j \rVert_{A_j}) +n\sum_{j=1}^{d}p_{\lambda_2}(w_{2j} \rVert b_j \rVert_{D_j}).
\end{equation}
To find the minimum of (\ref{q6}) for fixed tuning parameters, we use the iterative local quadratic approximation (LQA) proposed by Fan and Li (2001). Using a simple Taylor expansion, given an initial estimate $ b_{j}^{0} $, if $ \rVert b_j \rVert_{A_j}>0 ~~ \text{and}~~ \rVert b_j \rVert_{D_j}>0 $,
we approximate the penalty terms by
\begin{equation*}
p_{\lambda_1}(w_{1j} \rVert b_j \rVert_{A_j})\approx p_{\lambda_1}\left(w_{1j} \rVert b_{j}^{(0)} \rVert_{A_j}\right)+
\frac{1}{2}\frac{p'_{\lambda_1}\left( w_{1j}\rVert b_{j}^{(0)} \rVert_{A_j}\right)}{w_{1j} \rVert b_{j}^{(0)} \rVert_{A_j}}
\left\lbrace w_{1j}^{2}\rVert b_{j} \rVert_{A_j}^{2}-w_{1j}^{2}\rVert b_{j}^{(0)} \rVert_{A_j}^{2}\right\rbrace,
\end{equation*}
and
\begin{equation*}
p_{\lambda_2}\left(w_{2j} \rVert b_j \rVert_{D_j}\right)\approx p_{\lambda_2}\left(w_{2j} \rVert b_{j}^{(0)} \rVert_{D_j}\right)+
\frac{1}{2}\frac{p'_{\lambda_2}\left( w_{2j}\rVert b_{j}^{(0)} \rVert_{D_j}\right)}{w_{2j} \rVert b_{j}^{(0)} \rVert_{D_j}}
\left\lbrace w_{2j}^{2}\rVert b_{j} \rVert_{D_j}^{2}-w_{2j}^{2}\rVert b_{j}^{(0)} \rVert_{D_j}^{2}\right\rbrace.
\end{equation*}
After removing some irrelevant terms, the criterion becomes
\begin{equation}\label{q7}
Q(b)=\frac{1}{n} \rVert Y-Zb \rVert ^2 +\frac{1}{2}b^{T}(\Omega_1+\Omega_2)b
\end{equation}
for two $ dK \times dK $ matrices $ \Omega_1 $ and $ \Omega_2 $ defined by
\begin{equation}
\Omega_1=\text{diag}\left( \frac{\lambda_{1}w_{11}}{\rVert b_{1}^{(0)} \rVert_{A_1}}A_1,\dots,
\frac{\lambda_{1}w_{1d}}{\rVert b_{d}^{(0)} \rVert_{A_d}}A_d \right) \nonumber
\end{equation}
and
\begin{equation}
\Omega_2=\text{diag}\left( \frac{\lambda_{2}w_{21}}{\rVert b_{1}^{(0)} \rVert_{D_1}}D_1,\dots,
\frac{\lambda_{1}w_{2d}}{\rVert b_{d}^{(0)} \rVert_{D_d}}D_d \right) \nonumber
\end{equation}

Note that (\ref{q7}) is a quadratic function and thus there exists a
closed-form solution. Then the updating equation given the current estimate
$ b^{(0)} $ is
\begin{equation}\label{key}
b=\big(Z^{T}Z+n(\Omega_1+\Omega_2 )\big)^{-1}Z^{T}Y
\end{equation}

The algorithm repeatedly solves the minimization criterion (\ref{q7}) and updates $ b^{(m)} $ to $ b^{(m+1)} $, $ m=0,1,... $ until convergence. That is, in the m-th iteration, we solve (\ref{q7}), where $ \Omega_1 $ and $ \Omega_2 $ are as defined above but with $ b_{j}^{0} $ replaced by the current estimate $ b_{j}^{(m)} $. The solution obtained from (\ref{q7}) is the new estimate $ b^{(m+1)} $. \
During the iterations, as soon as some $ \rVert b_j \rVert_{A_j} $
(respectively, $ \rVert b_j \rVert_{D_j} $) drops below a certain threshold ($ 10^{-6} $ in our implementation), the component is identified as a zero function (respectively, linear function).\

\section{Simulation Studies}
For brevity, we refer to our approach as nonparametric correlation rank screening (NCRS). In this section, four simulation examples including different additive models with various scenarios are presented. The first three examples are allocated to our proposed screening procedure, while in the fourth one, the capability of structure identification method of Lian et al. (2015) is also examined.
In the former cases, the finite sample performance of the NCRS is compared with the existing competitors, such as the SIRS (Zhu et al., 2011), SIS (Fan and Lv, 2008), NIS (Fan et al., 2011) and the CR-SIS (Zhang et al., 2017). We consider two criteria for evaluating the performance as described in Zhu et al. (2011). The first criterion is the minimum model size (denoted by \textit{M}), that is the smallest number of covariates needed to ensure that all the active variables are selected. To get better inference, the quantiles 5\%, 25\%, 50\%, 75\% and 95\% quantiles of \textit{M} out of 200 replications were also presented. The second criterion is the proportion (denoted by \textit{S}) of truly active predictors that are identified by the screening procedure for a given model size in 200 replications, when the threshold $ \nu_n=\big[n/log(n)\big] $ is adopted. Note that the first criterion does not need to specify a threshold.  The more reliable screening procedure, the closer \textit{M} value to the number of active predictor and also the closer \textit{S} value to 1. \

We also conduct some Monte Carlo studies to assess the effectiveness of our two stage proposed method to separation of the linear and nonlinear components and to identify insignificant covariates simultaneously in partial linear additive models of non-polynomial (NP) dimensionality based on double penalization.\

To implement the procedures described in this paper, we need to find a data-driven procedure to choose the regularization parameters $ \lambda_1~\text{and}~ \lambda_2 $, and numbers of spline bases $ K $ . However, choosing different node sequence for different coefficients are computationally hard. To ease the computational burden, we fix $  K = 6 $  following Huang et al. (2010) and Lian et al. (2015).
To select the regularization parameters $ \lambda_1~\text{and}~ \lambda_2 $
simultaneously, we use the extended Bayesian information criterion (eBIC) of Chen and Chen (2008) that was developed for parametric models.
In our context, a natural eBIC-type criterion is defined by
\begin{equation}\label{q9}
\text{log}(\frac{1}{n} \rVert Y-Z\hat{b}_\lambda \rVert ^2)
+d_1 \frac{\text{log}(n/K)}{n/K}+d_2 \frac{\text{log}n}{n}+
\frac{d_1K+d_2}{n}\text{log} d,
\end{equation}
where $ \hat{b}_\lambda $ is the minimizer of (\ref{q6}) for given
$ \lambda=(\lambda_1,\lambda_2) $, $ d_1 $ is the number of components
estimated as nonparametric and $ d_2 $ is the number of
components estimated as parametric, both for the given $ \lambda $.
\vspace{.5cm}

\textbf{Example 1.} In the first example, we consider a classical linear model with varying squared multiple correlation coefficient $ R^2 $ and error distribution:
\begin{equation}
Y = c\beta^{T}X + \sigma\varepsilon,
\end{equation}
where $ \beta=(1, 0.8, 0.6, 0.4, 0.2, 0,\dots, 0)^T $ takes grid values, i.e., only the first five predictors are active. This example is adapted from Zhu et al. (2011). The ultrahigh-dimensional covariate $ X=(X_1,\dots,X_p) $ follows a multivariate normal distribution with mean $ 0 $ and the covariance matrix $ \Sigma=(\sigma_{ij})_{p\times p} $ with $ \sigma_{ii}=1 $ and $ \sigma_{ij} = 0.8^{|i-j|} $ for $ i\neq j $. We set $ \sigma^2=6.83 $ and considered two error $ \varepsilon $ distributions, a standard normal $ N(0, 1) $ and a t-distribution with 5 degree of freedom that has a heavy tail. We varied the constant $ c $ in front of $ \beta^{T}X $ to control the signal-to-noise ratio.
We choose $ c=0.5, 1~\text{and}~2 $, with the corresponding $ R^2 = 20\%, 50\% ~\text{and}~ 80\%. $ The sample size and the number of covariates are set to $ n=200, p=2000 $, respectively. For each scenario, based on $ 200 $ simulation runs, the results are given in Table \ref{tab1}. Each scenario is designed as a combination  of \textit{i)} distribution of errors, \textit{ii)} c-values \textit{iii)} screening method.\

From Table \ref{tab1}, when $ c=1 $ and $ 2 $, for both cases $ N(0,1) $ and $ t(5) $, all five screening methods perform equally well in most cases. In these setting, the corresponding $ M-\text{values}=5 $ shows that at least in 190 runs (quantile 95\%) out of $ 200 $, the five active covariates are appeared in the first $ 5 $ position of sorted lists. \

The difference between these methods is emerged when $ c=0.5 $. Although in this case, our proposed NCRS method, with $ S=0.92 $ and $ 0.70 $, is less accurate than SIS, SIRS and CR-SIS, but it is comparable to the others for the normal error. For the scenario including normal error and $ c=0.5 $, It is also worth noting that SIS (with $ M_{0.95} = 19 $) performs better than other methods. This is due to the fact that the true model is linear and the covariates are jointly normally distributed, which implies that the marginal projection is linear as well. However, for the heavy-tailed error, the performances of the CR-SIS and SIRS procedures are comparable. In this case, the NCRS method performs better than NIS method, particularly for $ c=0.5 $.
\vspace{.5cm}

\begin{table}[H]
	\caption{Five quantiles of minimum model size and the proportion of \textit{S}  among $ 200 $ replications in Example 1 with the true model size $ 5 $ and $ p=2000 $.}\label{tab1}
	\centering
	{\small \begin{tabular}{lllllllll}
			\hline	\hline
			$ \varepsilon $ & c & method &  & & M & &   & S\\
			\cline{4-8}
			& &  & 5\% & 25\% & 50\%  & 75\% & 95\%  & \\
			\hline	\hline   
 $ N(0,1) $ &0.5&NCRS & 5& 5 & 5 & 8 & 84 & 0.92\\
		&	&SIS (Fan et al. 2008)  & 5& 5 & 5 & 6 & 19 & 0.95\\
		&	&NIS (Fan et al. 2011)  & 5& 5 & 5 & 8 & 83 & 0.92\\
		&	&SIRS (Zhu et al. 2011)  & 5& 5 & 5 & 6 & 36 & 0.95\\
		&	&CR-SIS (Zhang et al. 2017)  & 5& 5 & 5 & 6 & 28 & 0.96\\
		
		& 1 &NCRS & 5& 5 & 5 & 5 & 5 & 1.00\\
		&	&SIS (Fan et al. 2008)  & 5& 5 & 5 & 5 & 5 & 1.00\\
		&	&NIS (Fan et al. 2011)  & 5& 5 & 5 & 5 & 5 & 1.00\\
		&	&SIRS (Zhu et al. 2011)  & 5& 5 & 5 & 5 & 5 & 1.00\\
		&	&CR-SIS (Zhang et al. 2017)  & 5& 5 & 5 & 5 & 5 & 1.00\\
		
		& 2 &NCRS & 5& 5 & 5 & 5 & 5 & 1.00\\
		&	&SIS (Fan et al. 2008)  & 5& 5 & 5 & 5 & 5 & 1.00\\
		&	&NIS (Fan et al. 2011)  & 5& 5 & 5 & 5 & 5 & 1.00\\
		&	&SIRS (Zhu et al. 2011)  & 5& 5 & 5 & 5 & 5 & 1.00\\
		&	&CR-SIS (Zhang et al. 2017)  & 5& 5 & 5 & 5 & 5 & 1.00\\	
			\hline 	
 $ t_{(5)} $ &0.5&NCRS & 5& 7 & 15 & 47 & 263 & 0.70\\
        &	&SIS (Fan et al. 2008) & 5& 6 & 9 & 20 & 212 & 0.82\\
        &	&NIS (Fan et al. 2011) & 5& 10 & 24 & 84 & 420 & 0.58\\
        &	&SIRS (Zhu et al. 2011) & 5& 5 & 6 & 11 & 104 & 0.90\\
        &	&CR-SIS (Zhang et al. 2017) & 5& 5 & 6 & 10 & 86 & 0.91 \\
   
        & 1 &NCRS & 5& 5 & 5 & 5 & 6 & 1.00\\
        &	&SIS (Fan et al. 2008) & 5& 5 & 5 & 5 & 5 & 1.00\\
        &	&NIS (Fan et al. 2011) & 5& 5 & 5 & 5 & 9 & 0.99\\
        &	&SIRS (Zhu et al. 2011) & 5& 5 & 5 & 5 & 6 & 1.00\\
        &	&CR-SIS (Zhang et al. 2017) & 5& 5 & 5 & 5 & 5 & 1.00\\
   
        & 2 &NCRS & 5& 5 & 5 & 5 & 5 & 1.00\\
        &	&SIS (Fan et al. 2008) & 5& 5 & 5 & 5 & 5 & 1.00\\
        &	&NIS (Fan et al. 2011) & 5& 5 & 5 & 5 & 5 & 1.00\\
        &	&SIRS (Zhu et al. 2011) & 5& 5 & 5 & 5 & 5 & 1.00\\
        &	&CR-SIS (Zhang et al. 2017) & 5& 5 & 5 & 5 & 5 & 1.00\\
			\hline	\hline
	\end{tabular} } 	
\end{table}


\textbf{Example 2.} Following Fan, Feng and Song (2011), we generate the data from the following additive model:
\begin{equation*}
Y=5g_{1}(X_{1})+3g_{2}(X_{2})+4g_{3}(X_{3})+6g_{4}(X_{4})+\sqrt{1.74} \varepsilon, 
\end{equation*}
where $ g_1(x)=x,~ g_2(x)=(2x-1)^2,~ g_3(x)=\sin(2\pi x)/(2-\sin(2\pi x)),~ g_4(x)=0.1 \sin(2\pi x)+0.2 \cos(2\pi x)+0.3 \sin(2\pi x)^2+0.4 \cos(2\pi x)^3+0.5 \sin(2\pi x)^3 $ and the vector of covariates $ X = (X_1,\dots, X_p)^T $ is generated in the same way as that in Example 1. We presented the simulation results for $ M $ and $ S $ in Table \ref{tab2}. 

\begin{table}[t]
	\caption{Five quantiles of minimum model size and the proportion of \textit{S} among 200 replications in Example 2 with the true model size $ 4 $ and $ p=2000 $.}\label{tab2}
	\centering
	{\small \begin{tabular}{lllllllll}
			\hline	\hline
		$ \varepsilon $ &	n & method &  & & M & &   & S\\
			\cline{4-8}
			&    &    & 5\% & 25\% & 50\%  & 75\% & 95\%   & \\
			\hline	\hline   
$ N(0,1)$ &200&NCRS & 4 & 4 & 4 & 4 & 5 & 1.00\\
		  & &NIS (Fan et al. 2011) & 4  & 4  & 4  & 4   & 5 &  1.00\\
	   	  & &SIRS (Zhu et al. 2011)& 4  & 8  & 29 & 132 & 681 &  0.53 \\
	      & &CR-SIS (Zhang et al. 2017) & 5  & 26 & 148 & 680 & 1476& 0.28\\

    	& 400 &NCRS & 4 & 4 & 4 & 4 & 4 & 1.00\\
		  & &NIS (Fan et al. 2011) & 4  & 4  & 4  & 4 & 4  &  1.00\\
	      &	&SIRS (Zhu et al. 2011)& 4  & 4  & 6 & 19 & 67 &  0.94 \\
	      &	&CR-SIS (Zhang et al. 2017) & 4  & 10 & 69 & 267 & 1187 & 0.49 \\
	      
	      \hline
$t_{(1)}$ &200&NCRS & 4 & 4 & 4 & 16 & 87 & 0.81\\
	      & &NIS (Fan et al. 2011) & 4  & 4  & 14  & 252 & 1796 &  0.56\\
	      & &SIRS (Zhu et al. 2011)& 4  & 10  & 43 & 184 & 883 &  0.48 \\
	      & &CR-SIS (Zhang et al. 2017) & 6 & 25 & 128 & 676 & 1563 &0.33\\
	      	
	      & 400 &NCRS & 4 & 4 & 4 & 4 & 4 & 1.00\\
	      & &NIS (Fan et al. 2011) & 4  & 4  & 4  & 4 & 4  &  1.00\\
	      &	&SIRS (Zhu et al. 2011)& 4  & 4  & 6 & 18 & 67 &  0.94 \\
	      &	&CR-SIS (Zhang et al. 2017) & 5  & 11 & 68 & 305 & 1124& 0.48\\
	      
	      \hline
$t_{(5)}$ & 200 &NCRS & 4 & 4 & 4 & 4 & 5  & 1.00\\
          & &NIS (Fan et al. 2011) & 4 & 4 & 4 & 4 & 5 & 1.00\\
          & &SIRS (Zhu et al. 2011)& 4 & 9 & 29 & 132 & 614  & 0.52\\
          & &CR-SIS (Zhang et al. 2017) & 5 & 28 & 126 & 659 & 1485 & 0.30\\
	      
	      & 400 &NCRS & 4 & 4 & 4 & 4 &4  & 1.00\\
	      & &NIS (Fan et al. 2011)& 4 & 4 & 4 & 4 &4  & 1.00\\
	      &	&SIRS (Zhu et al. 2011)&  4 & 4  & 6 & 19 & 78 &  0.94 \\
	      &	&CR-SIS (Zhang et al. 2017) & 4 & 11 & 61 & 270 & 1162 & 0.51\\
		\hline	\hline
	\end{tabular} } 	
\end{table}

From Table \ref{tab2}, for $ n=400 $ and for both types of distribution error, the action of the proposed NCRS and NIS, with $ M_{0.95}=4 $ and $ s=1.0 $, are very well and not comparable with the two others scenarios. As before, the $ M-\text{value}=4 $ implies that even in 95\% of time (among $ 200 $ runs), NCRS and NIS perfectly distinguish the four active covariates in the first four place of the sorted list. When $ n=200 $, these two methods provide nearly the same powerful results for normal and $ t(5) $ errors, but for the heavy tailed errors $ t(1) $, our NCRS is much superior than the others. In this setting, the significant difference between ``87", ``1796", ``883" and ``1563"  in the last column of $ M $ (95 \%) is an evidence for the superiority. Moreover, the corresponding $ S-\text{values}=0.81 $ for NCRS is much bigger than the others, which also shows the capability of our method. Both the above mentioned methods outperform SIRS and CR-SIS in must scenarios. 
\vspace{.5cm}

\textbf{Example 3.} This example is a more difficult case than Example's 1 and 2, because it has 8 important variables with different coefficients:
\begin{align*}
Y&=g_{1}(X_{1})+g_{2}(X_{2})+1.5g_{3}(X_{3})+1.5g_{4}(X_{4})+2g_{1}(X_{5})+2g_{2}(X_{6})\\
&+ 2.5g_{3}(X_{7})+2.5g_{4}(X_{8})+\varepsilon, 
\end{align*}
where $ g_j(x) $'s are the same as those in example 2. The ultrahigh-dimensional covariate $ X=(X_1,\dots,X_p) $ is generated from a multivariate normal distribution with mean $ 0 $ and the covariance matrix $ \Sigma=(\sigma_{ij})_{p\times p} $ with $ \sigma_{ii}=1 $ and $ \sigma_{ij} = 0.5 $ if both $ i,j\in \mathcal{A} $ or  $ i,j\in \mathcal{I} $,  and $ \sigma_{ij}=0.1 $ otherwise, where $ \mathcal{A} $  and $ \mathcal{I} $ are the active and inactive covariate sets, respectively. We presented the simulation results for $ M $ and $ S $ in Table \ref{tab3}.

According to Table \ref{tab3}, for the setting $ t(1) $ and $ t(5) $ distributions, by considering $ n=200 $, our NCRS is superior in terms of either the minimum model size required to cover all the active covariates or the proportion that all active predictors are selected. For the same setting and $ n=400 $, the NCRS is powerful as other competitors. In the case of normal errors, both NCRS and NIS have also similar performance and are equally well. Both of them outperform the SIRS and CR-SIS procedures.
\vspace{.5cm}

\begin{table}[t]
	\caption{Five quantiles of minimum model size and the proportion of \textit{S} among $ 200  $ replications in Example 3 with the true model size $ 8 $ and $ p=2000 $}\label{tab3}
	\centering
	{\small \begin{tabular}{lllllllll}
			\hline	\hline
			$ \varepsilon $ &	n & method &  & & M & &   & S\\
			\cline{4-8}
			&    &    & 5\% & 25\% & 50\%  & 75\% & 95\%   & \\
			\hline	\hline   
  $ N(0,1)$ &200&NCRS & 8 & 9 & 19 & 70 & 414 & 0.61\\
			& &NIS (Fan et al. 2011) & 8  & 8  & 13  & 48  & 380 &  0.71\\
			& &SIRS (Zhu et al. 2011)& 8 & 18  & 95 & 472 & 1219 &  0.33 \\
			& &CR-SIS (Zhang et al. 2016) & 10 & 112 & 421 & 1156 & 1913 & 0.14\\
			
			& 400 &NCRS & 8 & 8 & 8 & 8 & 13 & 1.00\\
			& &NIS (Fan et al. 2011) & 8  & 8  & 8  & 8 & 12  &  1.00\\
			& &SIRS (Zhu et al. 2011)& 8  & 8  & 10 & 26 & 266 &  0.81 \\
			& &CR-SIS (Zhang et al. 2017) & 8 & 30 & 167 & 719 & 1612 & 0.34\\
			
			\hline	
  $t_{(1)}$ &200&NCRS & 9 & 22 & 130 & 900 & 1834  & 0.33\\
			& &NIS (Fan et al. 2011) & 8  & 42 & 271 & 1618 & 1958 &  0.24\\
			& &SIRS (Zhu et al. 2011)& 8  & 42 & 199 & 626 & 1462 &  0.25 \\
			& &CR-SIS (Zhang et al. 2016)& 14 & 81 & 544 & 1232 & 1873 & 0.17 \\
			
		  & 400 &NCRS & 8 & 8 & 10 & 128 & 1648 & 0.70\\
			&   &NIS (Fan et al. 2011) & 8  & 8  & 31  & 628 & 1879 & 0.55\\
			&	&SIRS (Zhu et al. 2011)& 8  & 8  & 12 & 64 & 425 & 0.76 \\
			&	&CR-SIS (Zhang et al. 2017) & 8  & 17 & 113 & 638 & 1637 & 0.44\\
			
			\hline
 $t_{(5 )}$ &200&NCRS & 8 & 9 & 16 & 74 & 547 & 0.63\\
			& &NIS (Fan et al. 2011) & 8  & 8 & 11 & 46 & 479 &  0.73\\
			& &SIRS (Zhu et al. 2011)& 8  & 20 & 93 & 363 & 1305 &  0.31 \\
			& &CR-SIS (Zhang et al. 2016)& 9 & 110 & 491 & 1225 & 1760 & 0.14 \\
			
			& 400 &NCRS & 8 & 8 & 8 & 8 & 19 & 0.97\\
			&   &NIS (Fan et al. 2011) & 8  & 8  & 8  & 8 & 10 & 0.98\\
			&	&SIRS (Zhu et al. 2011)& 8  & 8  & 10 & 34 & 248 & 0.81 \\
			&	&CR-SIS (Zhang et al. 2017) & 8  & 17 & 137 & 638 & 1740 & 0.39\\
			\hline	\hline
	\end{tabular} } 	
\end{table}

\textbf{Example 4.} In this example, we first apply the NCRS method to reduce dimensionality, and then fit two models, a sparse additive model (SAM) where only one penalty is used to identify nonzero components (and thus parametric components cannot be identified), and a partial linear model (PLAM) where two penalty is used to simultaneously identify nonzero and linear components (Lian, et.al, 2015). 
We generated data from the model 
 \begin{equation}
Y=\sum_{j=1}^{p} f_{j}(X_{j})+\varepsilon,
\end{equation}
where $f_1(x)=5sin(2\pi x),~f_2(x)=10x(1-x),~ f_3(x)=3x,~f_4(x)=2x,~ f_5(x)=-2x, f_j(x)=0, j>5 $.
To generate covariates, we first let $  X_j $ be marginally standard normal with correlations
given by $ Cov(X_i,X_j) =0.8^{|i-j|} $, and then apply the cumulative distribution function of the standard normal
distribution to transform $ X_j $ to be marginally uniform on $ [0,1] $.
The noises are generated from mean zero normal distribution with standard deviation $ \sigma $. We performed simulations with $ n=200, 400, p=1000, 2000, $ and $ \sigma = 0.2, 0.5 $, resulting in eight scenarios.\

For all scenarios, 200 datasets are generated and the results are summarized in Table \ref{tab4}. 
We used several criterion to measure the model identification performance:
``NV": average number of variables selected; ``NVT":average number of variables selected that are truly significant; ``NN": average number of nonlinear components selected; ``NNT":average number of nonlinear components selected that are truly nonlinear; NL: average number of linear components selected; ``NLT":average number of linear components selected that are truly linear. The true number of nonparametric components is 2 and the true number of linear components is 3. The numbers in parenthesis are the corresponding standard errors. In terms of identifying the significant variables, the two methods perform similarly. However, the SAM cannot detect the parametric components.
\begin{table}
	\caption{Model identification results for Example 4. }\label{tab4}
	\centering
	{\small \begin{tabular}{lllllllll}
			\hline	\hline
		      & & SAM & & &  &  PLAM &   &  \\
		\cline{3-4} \cline{6-9}
            $ p $ &	$ n,\sigma $ & NV & NVT && NN & NNT & NL &  NLT\\
			\hline	\hline   
     1000 & $ n=200,\sigma=0.2 $ & 5.02(0.22) &4.98(0.12) &       & 2.33(0.77) & 1.94(0.42) & 2.75(0.87) & 2.59(0.69) \\
		   	&$ n=200,\sigma=0.5 $& 5.33(0.85) &4.90(0.32) &       & 2.34(1.02) & 1.80(0.60)  & 2.92(1.71) & 2.36(0.79) \\
			&$ n=400,\sigma=0.2 $& 5(0)       & 5(0)      &       & 2.06(0.64) & 1.98(0.36) & 2.96(0.74) & 2.92(0.48)\\
			&$ n=400,\sigma=0.5 $& 5.01(0.07) & 5(0)      &       & 2.14(0.79) & 1.86(0.52) & 2.87(0.78) & 2.68(0.58) \\
			
	\hline
    2000 & $ n=200,\sigma=0.2 $ & 4.97(0.20) &4.96(0.18) &       & 2.39(0.84) & 1.93(0.37) & 2.61(0.88) & 2.50(0.76) \\
	       &$ n=200,\sigma=0.5$ & 5.08(0.32) &4.98(0.12) &       & 2.40(1.21) & 1.75(0.66) & 2.85(1.83) & 2.23(0.88) \\
	       &$ n=400,\sigma=0.2$ & 5(0)       & 5(0)      &       & 1.94(0.69) & 1.83(0.54) & 3.01(0.76) & 2.86(0.53) \\
	       &$ n=400,\sigma=0.5$ & 5.02(0.08) & 5(0)      &       & 1.94(0.92) & 1.71(0.72) & 2.88(1.01) & 2.60(0.77) \\
	 		
	\hline	\hline
	\end{tabular} } 	
\end{table}

\section{Cardiomyopathy data}
In this section, we apply and evaluate our method to do variable selection and identify the structure of components, for a real dataset. This dataset was analyzed by Segal et al. (2003), Hall and Miller (2009).  The aim is to identify the most influential genes for overexpression of a G protein-coupled receptor, designated Ro1, in mice. The Ro1 expression level, $ Y_i $, was measured for $ n = 30 $ specimens, and genetic expression levels, $ X_i $, were obtained for $  p = 6,319 $ genes.\

\begin{figure}[H]
	\begin{center}
		\includegraphics[width=13cm,height=7cm]{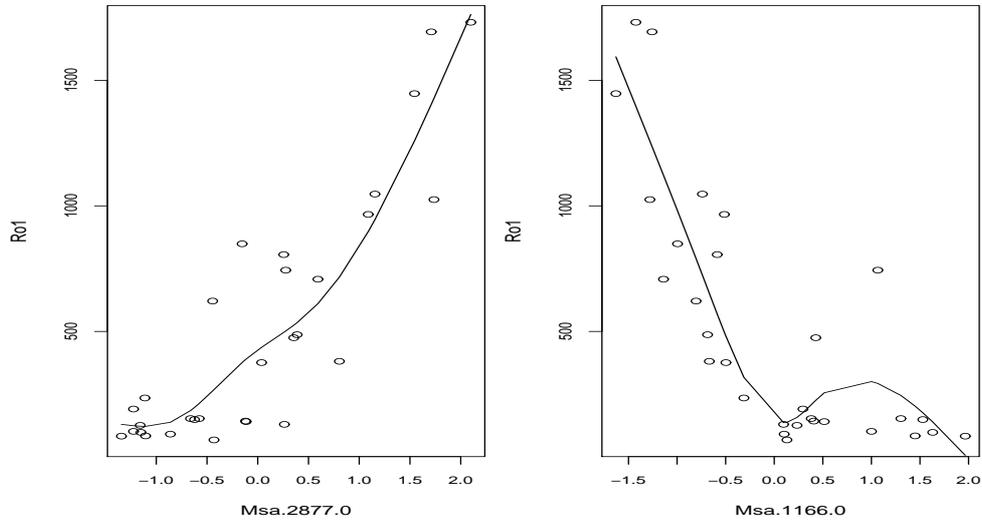}
		\caption{\footnotesize{The scatter plot of Y versus two gene expression levels with cubic spline fit curves.}}\label{fig1}
	\end{center}
\end{figure}

According to Figure \ref{fig1}, the scatterplots of $ Y $ versus these two gene expression levels with cubic spline fit curves indicate clearly the existence of nonlinear patterns. \
The NCRS procedure ranks two genes, Msa.2134.0 and Msa.2877.0, in the top, which is the same as CR-SIS (Zhang et al., 2017) and SIRS (Zhu et al.,2011) do. The NIS procedure ranks two genes, labeled as  Msa.2877.0 and Msa.1166.0, at the top.\

\begin{figure}[H]
	\begin{center}
		\includegraphics[width=14cm,height=11cm]{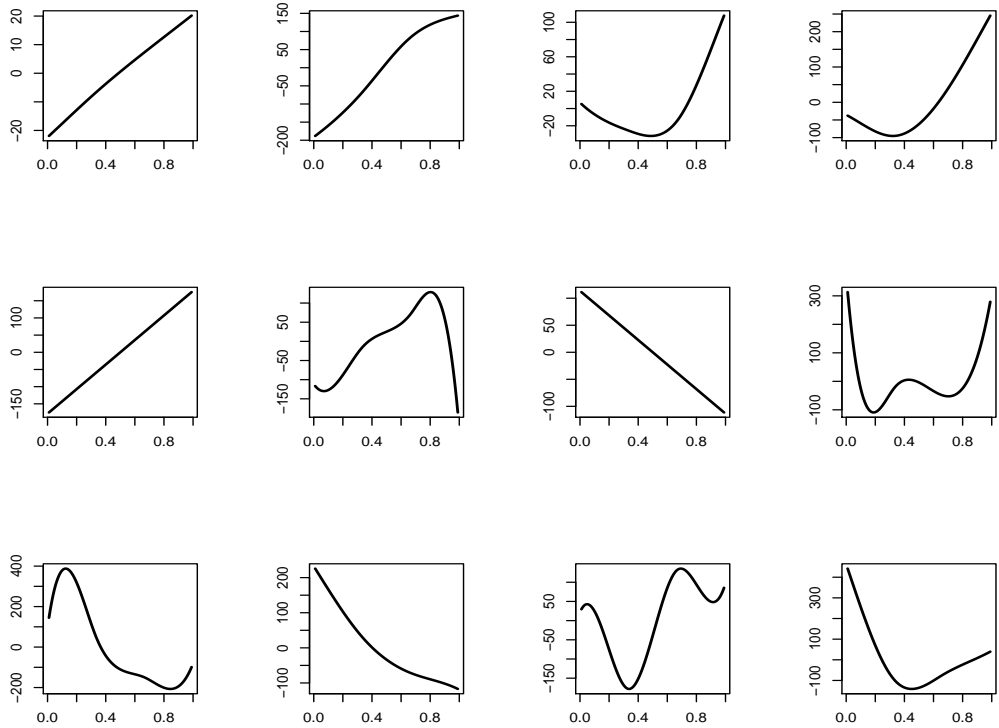}
		\caption{\footnotesize{Fitted regression functions for the 12 genes selected.}}\label{fig2}
	\end{center}
\end{figure}

We first applied NCRS to reduce the covariate dimension to the size
 $ 2[n/log(n)]=16 $, and then obtained the sparse additive estimator and partial linear additive estimator. Both of them selected the significant variables similarly. After applying double penalization based  procedure, we have identified 3 genes of linear effects and 9 genes of nonlinear effects. The genes of linear effects are Msa.10108.0, Msa.2134.0, and Msa.26025.0, whereas, the genes of nonlinear effects are Msa.1166.0 and Msa.15405.0, Msa.1590.0, Msa.2400.0, Msa.2877.0, Msa.5583.0, Msa.5794.0, Msa.7336.0. Their effect functions are depicted in Figure \ref{fig2}.

 To evaluate the performance of the sparse additive and partial linear additive models, we used leave one out cross validation and compared the prediction mean squared errors (PE). The PE for sparse additive model is 0.86 and for partial linear additive model is 0.83. Apparently the partial linear additive model has a smaller PE indicating a satisfactory prediction performance. 

\section*{Discussion}
In this article, we proposed a sure independence screening procedure in partially linear additive models using covariance between marginal nonparametric functions and the unconditional distribution function of  
$ Y $, that is, NCRS. We used B-spline basis functions for fitting the marginal nonparametric components. We also established the sure screening property for this procedure under some conditions. Moreover, in order to distinguish linear and nonlinear parts and to identify insignificant covariates simultaneously, we used a double penalization based procedure.
We examined the finite sample performance of the proposed procedure via an extensive Monte Carlo study and evaluated the proposed methodology through the analysis of Cardiomyopathy microarray dataset. Numerical studies demonstrated a satisfactory performance of our screening procedure and it is competitive with the existing procedures such as the SIS, NIS, SIRS and CR-SIS procedures.\

Similar to the SIS, the NCRS may fail to identify some important predictors that are jointly but not marginally important. Thus, it is of interest to develop an iterative procedure to fix such an issue. We used only the adaptive lasso penalty but other penalties such as smoothly clipped absolute deviation (SCAD, Fan and Li; 2001) and minimax concave penalty (MCP, Zhang; 2010) could also be applied. 

\section*{Appendix}

\textbf{Proof of Theorem 1:} Let

\begin{equation}
r_{k}^{*}=\Big\{\frac{1}{n}\sum_{i=1}^n m_{k}(X_{ik}) G(Y_i)\Big\}^2.
\end{equation} 
 We prove this theorem via two steps. First, we derive the exponential tail probability bound of $ p(|\hat{r}_k-r_{k}^{*}|\geq \nu n^{-\alpha}) $ for any positive constants $ \nu $ and $ 0\leqslant\alpha < 1/2 $. Straightforward calculations entail that
 \begin{align*}
\mid\hat{r}_k-r_{k}^{*}\mid&=\mid\Big\{\frac{1}{n}\sum_{i=1}^n \hat{m}_{k}(X_{ik}) \hat{G}_n(Y_i)\Big\}^2-\Big\{\frac{1}{n}\sum_{i=1}^n m_{k}(X_{ik}) G(Y_i)\Big\}^2\mid\\
 &=\mid\Big(\frac{1}{n}\sum_{i=1}^n \hat{m}_{k}(X_{ik}) \hat{G}_n(Y_i)+\frac{1}{n}\sum_{i=1}^n m_{k}(X_{ik}) G(Y_i)\Big)\\
 &\Big(\frac{1}{n}\sum_{i=1}^n \hat{m}_{k}(X_{ik}) \hat{G}_n(Y_i)-\frac{1}{n}\sum_{i=1}^n m_{k}(X_{ik}) G(Y_i)\Big)\mid.
 \end{align*}
 By the SLLN, we have $  \frac{1}{n}\sum_{i=1}^{n}m_k(X_{ik})^2\xrightarrow{a.s}E[m_k(X_{ik})]^2. $
 Combining it with condition $ C1 $, there exists a positive constant $ c_1 $ such that
 \begin{equation}\label{ap1}
 \frac{1}{n}\sum_{i=1}^{n}m_k(X_{ik})^2\leq c_{1}^2.
 \end{equation}
holds a.s. when $ n $ is sufficiently large. Without loss of generality, assume that (\ref{ap1}) holds for the total probability space as the set with measure zero does not affect the derivations. Using the Cauchy-Schwarz inequality and the boundedness of $ \hat{G}_n(t) $ and $ G(t) $, we have 
 \begin{equation}\label{ap2}
 \mid\frac{1}{n}\sum_{i=1}^n \hat{m}_{k}(X_{ik}) \hat{G}_n(Y_i)\mid\leq c_1
 \qquad\text{and}\qquad
 \mid\frac{1}{n}\sum_{i=1}^n m_{k}(X_{ik}) G(Y_i)\mid\leq c_1.
 \end{equation}
 Using (\ref{ap1}) and (\ref{ap2}), we have
 \begin{align*}
 \mid\hat{r}_k-r_{k}^{*}\mid &\leq c_3\mid\frac{1}{n}\sum_{i=1}^n \hat{m}_{k}(X_{ik}) \hat{G}_n(Y_i)-\frac{1}{n}\sum_{i=1}^n m_{k}(X_{ik}) G(Y_i)\mid\\
 &\leq c_3\mid\frac{1}{n}\sum_{i=1}^n \hat{m}_{k}(X_{ik})\Big( \hat{G}_n(Y_i)-G(Y_i)\Big)\mid+\mid\frac{1}{n}\sum_{i=1}^n G(Y_i)\Big( \hat{m}_{k}(X_{ik})-m_{k}(X_{ik})\Big)\mid\\
 &\leq c_3\Big(\frac{1}{n}\sum_{i=1}^n \hat{m}_{k}^2(X_{ik})\Big)^{\frac{1}{2}}\max_{1\leq i\leq n}\mid \hat{G}_n(Y_i)-G(Y_i)\mid+\Big(\frac{1}{n}\sum_{i=1}^n ( \hat{m}_{k}(X_{ik})-m_{k}(X_{ik}))^2\Big)^\frac{1}{2}\\
 &= c_3 O(1/2)\max_{1\leq i\leq n}\mid \hat{G}_n(Y_i)-G(Y_i)\mid+O(1/2)\\
 & \leq c_4\max_{1\leq i\leq n}\mid \hat{G}_n(Y_i)-G(Y_i)\mid+O(1/2)\\ 
 & \leq c_4\max_{y \in \mathbb{R}}\mid \hat{G}_n(y)-G(y)\mid+O(1/2),
 \end{align*}
 where
 $c_3=2c_1 ~\text{and}~ c_4=c_3 O(1/2) $. It follows from the Dvoretzky-Kiefer-Wolfowitz inequality that
 \begin{align}\label{ap3}
 p\Big(\mid\hat{r}_k-r_{k}^{*}\mid\geq \nu n^{-\alpha} \Big)
 &\leq p\Big(c_4\max_{y \in \mathbb{R}}\mid \hat{G}_n(y)-G(y)\mid\geq \nu n^{-\alpha}-O(1/2) \Big) \nonumber\\ 
 &\leq 2\exp\{-2nc_4^{-2}( \nu n^{-\alpha}-O(1/2) )^2\}.
 \end{align}
 
 Second, we derive the exponential tail probability bound of $ p\Big(\mid r_{k}^{*}-r_{k}\mid\geq \nu n^{-\alpha} \Big) $ for any positive constants $ \nu $ and $ 0 \leq\alpha < 1/2. $ Using the similar arguments, we also have
 \begin{equation*}
 \mid r_{k}^{*}-r_{k}\mid \leq c_3\mid\frac{1}{n}\sum_{i=1}^n m_{k}(X_{ik}) G_n(Y_i)-E\{m_{k}(X_{k}) G(Y)\}\mid.
 \end{equation*}
 By the exponential Chebyshev inequality, for any $\xi > 0 $, we have
 
 \begin{align}
 p\Big(\mid r_{k}^{*}-r_{k}\mid \geq \nu n^{-\alpha} \Big)
 &\leq p\Big(c_3\mid\frac{1}{n}\sum_{i=1}^n m_{k}(X_{ik}) G(Y_i)-E\{m_{k}(X_{k}) G(Y)\}\mid\geq \nu n^{-\alpha} \Big)\nonumber\\ 
 &=p\Big(\mid\frac{1}{n}\sum_{i=1}^n m_{k}(X_{ik}) G(Y_i)-E\{m_{k}(X_{k}) G(Y)\}\mid \geq c_3^{-1} \nu n^{-\alpha}\Big)\nonumber \\ 
 &\leq \exp(-\zeta c_3^{-1}\nu n^{-\alpha}).E\Big(
 \exp\{\zeta \mid\frac{1}{n}\sum_{i=1}^n m_{k}(X_{ik}) G(Y_i)-E\{m_{k}(X_{k}) G(Y)\}\mid\}\Big)
 \end{align}
 
 Using the law of the iterated logarithm, we have
 
 \begin{equation}\label{ap4}
 \limsup_{n\rightarrow \infty}\frac{\sum_{i=1}^n m_{k}(X_{ik})G(Y_i)-nE\{m_{k}(X_{k})G(Y)\}}{\Big[n\log\log n.Var\{m_{k}(X_{k})G(Y)\}\Big]^{\frac{1}{2}}}=\sqrt{2}, \qquad \text{a.s.}
 \end{equation}
 Without loss of generality, when $ n $ is large enough and removing a zero measure set, under condition $ C1 $, there exists a positive constant $  c_5 $ such that
 \begin{equation*} \label{ap5}
 \Big(\frac{n}{\log\log n}\Big)^{\frac{1}{2}}.\Big[\frac{1}{n}\sum_{i=1}^n m_{k}(X_{ik})G(Y_i)-E\{m_{k}(X_{k})G(Y)\}\Big]\leq c_5.
 \end{equation*}
 We chose $ \zeta = \Big(\frac{n}{\log\log n}\Big)^{\frac{1}{2}} $, then it follows from (\ref{ap4}) and (\ref{ap5}) that
 \begin{equation}\label{ap6}
 p\Big(\mid r_{k}^{*}-r_{k}\mid \geq \nu n^{-\alpha} \Big)\leq
 \exp\Big\{-\Big(\frac{n^{1-2\alpha}}{\log\log n}\Big)^{\frac{1}{2}}c_3^{-1}c_5\nu\Big\}
 \end{equation}
 Combining (\ref{ap3}) and (\ref{ap6}), we have
\begin{align}\label{ap7}
p\Big(\mid \hat{r}_{k}-r_{k} \mid \geq 2\nu n^{-\alpha} \Big)
&\leq p\Big(\mid \hat{r}_{k}-r_{k}^{*}\mid \geq \nu n^{-\alpha} \Big)+p\Big(\mid r_{k}^{*}-r_{k}\mid \geq \nu n^{-\alpha} \Big)\nonumber\\
&\leq 2\exp\Big\{-2nc_4^{-2}( \nu n^{-\alpha}-O(1/2) )^2\Big\}+\exp\Big\{-\Big(\frac{n^{1-2\alpha}}{\log\log n}\Big)^{\frac{1}{2}}c_3^{-1}c_5\nu\Big\}\nonumber\\
& \leq O\Big[\exp\Big\{-\eta\Big(\frac{n^{1-2\alpha}}{\log\log n}\Big)^{\frac{1}{2}}\Big\}\Big]
\end{align}
where $ \eta =c_{3}^{-1} c_{5}\nu $. Immediately, we have
\begin{equation}\label{ap8}
p\Big(\max_{1\leq k \leq p}\mid \hat{r}_{k}-r_{k} \mid \geq 2\nu n^{-\alpha} \Big) \leq O\Big[p\exp\Big\{-\eta\Big(\frac{n^{1-2\alpha}}{\log\log n}\Big)^{\frac{1}{2}}\Big\}\Big]
\end{equation}
which proves the first part of Theorem 1 by taking $ c=2\nu $.
If $ \mathcal{A}\varsubsetneq \mathcal{\hat{A}} $, then there must exist some $ k \in \mathcal{A} $ such that $ \hat{r}_k<cn^{-\alpha} $. 
It follows from condition C2 that $ \mid\hat{r}_k-r_k\mid>cn^{-\alpha} $ for some $ k \in \mathcal{A} $, which implies that
$ \{\mathcal{A}\subsetneq\mathcal{\hat{A}}\} \subseteq\{\mid\hat{r}_k-r_k\mid>cn^{-\alpha}~\text{for some}~ k\in\mathcal{A} \}$.
 As a result, $\{\max_{k \in \mathcal{A}}\mid\hat{r}_k-r_k\mid\leq cn^{-\alpha} \}\subseteq \{\mathcal{A}\subseteq \mathcal{\hat{A}}\} $. 
Using (\ref{ap7}), we have
\begin{equation*}
p(\mathcal{A}\subseteq\mathcal{\hat{A}})\geq p\Big(\max_{k \in \mathcal{A}}\mid\hat{r}_k-r_k\mid\leq cn^{-\alpha}\Big)\\
\geq1-O\Big[a_n\exp\Big\{-\eta\Big(\frac{n^{1-2\alpha}}{\log\log n}\Big)^{\frac{1}{2}}\Big\}\Big]
\end{equation*}
where $ a_n = |A|$. Thus, the proof of Theorem 1 is completed.

\end{document}